# Enseignement des mathématiques, pédagogie concrète et méthodes actives dans les centres d'apprentissage (1945-1953)

Xavier SIDO
Univ. Lille, EA 4354 – Théodile-CIREL – Centre interuniversitaire de Recherche en Éducation de Lille, F59000 Lille, France

C'est en 1945 avec la création des centres d'apprentissage que se met en place une formation professionnelle scolarisée des ouvriers et des employés qualifiés. Ces établissements préparent en trois ans au certificat d'aptitude professionnelle (CAP), qui certifie essentiellement la maîtrise de savoir-faire pratiques. Dans la hiérarchie des établissements de l'enseignement technique (ou simplement ET), fondée sur le niveau de qualification auquel conduisent les établissements techniques, les centres d'apprentissage occupent le premier échelon en formant l'enseignement technique court. Les collèges technique (CT) et les écoles nationales professionnelles (ENP), **qui forment en 6 ou 7 ans (??)** à des emplois nécessitant davantage de polyvalence technique et de connaissances théoriques que de savoir-faire pratiques (employés hautement qualifiés, techniciens), constituent l'enseignement technique long.

Les centres d'apprentissage s'inscrivent dans la perspective d'un enseignement post-primaire. Ils sont en effet destinés à des adolescents issus des classes de fin d'études primaires qui, sinon, ne poursuivraient pas leurs études et intègreraient le monde du travail sans réelle formation professionnelle. Afin d'adapter la formation à ce nouveau public scolaire, de nombreuses enquêtes et études sont menées à la Libération[1]. Elles indiquent que ce public est d'une origine sociale plus modeste que celui des autres établissements techniques. De plus, l'absence de sélection à l'entrée de la plupart des centres[2], à la différence des autres écoles techniques, conduit des élèves d'un faible niveau en mathématiques à intégrer ces établissements.

L'enseignement dans les centres, qui vise à former selon une formule de l'époque, l'homme, le travailleur et le citoyen, associe formation professionnelle et formation générale. L'enseignement de mathématiques alors mis en place s'inscrit dans ces missions. Ses finalités sont de participer à l'apprentissage d'un métier et à l'acquisition par les élèves d'une certaine culture. Par culture, il faut entendre les connaissances et la formation de l'esprit nécessaires pour que ces derniers puissent s'intégrer dans la société et dépasser leur condition de travailleur. Dans la perspective d'une transition entre l'école et la vie active, les programmes sont largement influencés par ceux des classes de fin d'études tout en intégrant des contenus spécifiques nécessaires pour bien suivre les disciplines professionnelles. Le public des élèves et les exigences certificatives, fondateurs et structurants pour l'enseignement de mathématiques dans les centres[3], conduisent ainsi à mettre en place un enseignement spécifique[4] par rapport aux enseignements de mathématiques en vigueur dans les autres filières scolaires. Proche du modèle de l'enseignement dispensé dans les écoles primaires élémentaires, il prépare comme lui ses élèves à la vie courante et professionnelle. Cet enseignement éminemment pratique et utilitaire s'en distingue toutefois par les liens qu'il entretient avec les enseignements professionnels.

---

[1] Paul JUIF, « Une enquête sur les Centres d'Apprentissage », *Apprentissage, 1*, 1945, p. 3-15. ; Paul JUIF, « Une enquête sur les Centres d'Apprentissage », *Apprentissage, 2*, 1945, p. 81-85.

[2] Patrice PELPEL et Vincent TROGER, *Histoire de l'enseignement technique*, Paris : L'harmattan, 2001, p. 172.

[3] Maryse LOPEZ et Xavier SIDO, « L'enseignement des mathématiques et du français dans l'enseignement technique court de 1945 à 1985. Identité singulière, dynamique et temporalité spécifiques ? » dans R. D'ENFERT et J. LEBEAUME (dir), *Réformer les disciplines scolaires*, Grenoble : PUG, 2015, p. 137-154.

[4] Xavier SIDO, « L'enseignement des mathématiques en CAP dans l'enseignement professionnel scolarisé, 1945-1985 » dans G. BRUCY, F. MAILLARD et G MOREAU (coord.), *Le CAP, 1911-2011. Formations professionnelles, certifications et sociétés,* Rennes : PUR, 2013, p. 89-102.

Enseignement des mathématiques, pédagogie concrète et méthodes actives
dans les centres d'apprentissage (1945-19453)

En 1945, les missions des centres, la forte présence de maîtres issus du primaire au moment de leur mise en place, leur recrutement ainsi que le flou institutionnel relatif à leur inscription dans le système éducatif[5] posent la question de leur inscription dans le système d'enseignement. Font-ils partie d'un primaire prolongé comme semble le suggérer l'inspection de l'enseignement primaire, ou d'un secondaire technique comme le revendiquent les cadres de l'ET ?[6] Les débats renvoient à des enjeux corporatistes, mais aussi pédagogiques. « Existe-t-il aussi des problèmes pédagogiques propres à notre enseignement général, problèmes qui ne se posent ni dans l'enseignement primaire, ni dans le reste du second degré ? »[7], demande Fernand Canonge, professeur en ENNA[8], en 1949. Ce questionnement intervient dans une période marquée à la fois par une intense réflexion sur les méthodes d'enseignement dans l'ET et par l'influence du mouvement des pédagogies nouvelles qui défend l'idée d'une participation active des individus à leur formation, dans une démarche dont le point de départ serait leurs centres d'intérêts[9].

Dans ce contexte, l'enjeu de cette contribution est d'identifier et de caractériser le modèle pédagogique prescrit pour cet enseignement mathématique particulier dispensé dans le cadre d'une formation post-primaire professionnelle, à savoir les centres d'apprentissage. Plus spécifiquement, nous visons à répondre à deux questions qui organisent notre texte : quelles spécificités, proximités, différences avec les modèles pédagogiques des enseignements mathématiques pour le primaire élémentaire et le secondaire ? Quelle opérationnalisation de cet enseignement mathématique, notamment au regard des spécificités du public auquel il est destiné ?

---

[5] Il faut attendre la Loi du 21 février 1949 relative au statut des centres d'apprentissage (JO n°46 du 22 février 1949, p. 1943-1944) pour que ces derniers soient pleinement intégrés à l'ET.

[6] Fernand, CANNONGE, « Les centres font partis du second degré », *Technique, Art, Science, 1,* 1946, p.7-9.

[7] Fernand CANNONGE, « Les Écoles Normales Nationales d'Apprentissage : les problèmes pédagogiques propres à l'Enseignement Technique », *Technique, Art, Science, 3,* 1949, p. 38-42, p. 38.

[8] École normale nationale d'apprentissage. Mises en place à partir de 1946, les ENNA sont chargées de la formation des enseignants des centres d'apprentissage.

[9] Laurent BESSE, Laurent GUTIERREZ et Antoine PROST (dir.), réformer l'école : l'apport de l'École, nouvelle (1930-1970), Grenoble : PUG, 2012.

Les sources de l'enquête sont les discours pédagogiques relatifs à l'enseignement mathématique dans les centres et ceux plus généraux fixant les orientations de l'ET. Ils peuvent être appréhendés dans les premières instructions officielles pour les centres (1945) ainsi que dans les articles publiés dans les revues para-officielles *Technique, Art, Science* (*TAS*) et *Apprentissage*, crées respectivement en 1946 et 1945. Ces deux revues publient des articles chargés de guider les professeurs des centres qui ont valeur de directives officielles. Ce sont, par exemple, des articles de pédagogie ou des leçons-types de mathématiques dont la plupart concernent les centres et qui sont écrits principalement par des professeurs d'ENNA ou de centre d'apprentissage et par des inspecteurs de l'ET.

Ces sources délimitent notre période d'étude, qui va de 1945 à 1953. Si 1945 correspond à la création des centres d'apprentissage, 1953 marque la parution d'un programme définitif pour cet enseignement qui remplace ainsi les premières instructions provisoires de 1945. De plus, cette date marque le tarissement, dans les deux revues mentionnées plus haut, des articles relatifs à l'enseignement de mathématiques dans les centres. Les raisons en sont sans doute la systématisation de la formation des enseignants nouvellement recrutés dans les ENNA à partir de 1950 et le départ du responsable de la section mathématiques dans la revue *TAS*, René Cluzel, auteur de nombreux manuels pour l'enseignement de mathématiques dans les centres.

## Une démarche pédagogique spécifique en tension entre le modèle primaire et secondaire

Pour que l'enseignement mathématique dans les centres puisse contribuer à rendre accessible aux élèves une certaine forme de culture et participer à leur épanouissement, inspecteurs, professeurs d'ENNA et de centres pensent qu'il ne doit pas contrarier leur nature :

> *On n'y réussira pas, même avec les meilleures intentions de bien faire si l'on contrarie les tendances profondes de l'être. Il faut connaître son évolution mentale, agir à temps en suivant la pente de ses intérêts et la direction générale de sa personnalité* [10]

Si le constat des difficultés des élèves à acquérir les rudiments du raisonnement et à maîtriser les notions abstraites est unanimement partagé dans les revues, deux discours, parfois étroitement articulés, tendent à l'expliquer et à y remédier. Ils s'appuient respectivement sur les centres d'intérêts spécifiques des élèves et sur leur jeune âge.

---

[10] Fernand CANONGE, « Le point de vue de l'ET sur les problèmes de l'orientation scolaire », *Technique, Art, Science, 1,* 1947, p.19-21.

Le premier discours est porté principalement par des inspecteurs et certains maîtres des centres. Ils considèrent que si les élèves ont des difficultés en mathématiques, c'est qu'ils ne participent pas assez à la leçon à l'établissement des propriétés et des raisonnements. Plus précisément, ils considèrent que les raisons de leur « insuccès […] tiennent surtout à la peine de l'enfant à apprendre des définitions abstraites qu'il ne rattache à rien, à suivre des démonstrations souvent compliquées dont il ne voit pas la nécessité »[11]. Les critiques portent sur les méthodes « traditionnelles », considérées comme expositives, dogmatiques, fondées sur l'abstraction et qui seraient responsables d'une perte de sens dans les apprentissages mathématiques car inadaptées au profil des élèves. Ils fondent leur raisonnement sur l'impact qu'aurait l'apprentissage d'un métier sur ces adolescents qui se différencieraient alors de ceux qui fréquentent les lycées et les collèges par leur goût pour le concret et leur curiosité naturelle pour tout ce qui touche au métier[12]. Le second discours est tenu essentiellement par les professeurs d'ENNA et des enseignants des matières générales. Se référant à la psychologie de l'enfant et la psychopédagogie[13], ils réfutent le présupposé selon lequel « le raisonnement déductif correspond à une fonction primitive de la pensée des apprentis »[14], surtout au plus jeune âge. Ils militent en faveur de l'adaptation de l'enseignement aux différents stades de développement des élèves et donc préconisent pour les mathématiques de différencier les méthodes pédagogiques en fonction de l'âge mental des élèves.

Ainsi, en référence à la fois aux centres d'intérêt des élèves et à leur jeune âge, les instructions et les publications pédagogiques recommandent que le concret soit le point de départ de l'apprentissage des mathématiques. Afin d'expliciter ce modèle pédagogique, nous reprenons ici les trois phases que François Matray, professeur d'ENNA, distingue pour l'enseignement des mathématiques dans son ouvrage *Pédagogie de l'enseignement technique*, paru en 1952[15]. Elles rendent compte de la progression de l'enseignement dans les programmes durant les trois ans que durent les études dans les centres, mais aussi de la structure interne des leçons proposées dans les revues.

---

[11] LABARTHE, « L'enseignement de la Géométrie », *Technique, Art, Science, 7,* 1947, p. 44-45, p. 44.

[12] Fernand CANONGE, « Les bases psychologiques de notre enseignement », *Technique, Art, Science, 3,* 1948, p. 5-8.

[13] Respectivement les travaux de Jean Bourjade, spécialiste de la psychologie de l'enfant et directeur de l'école pratique et de pédagogie de Lyon et Louis Dintzer professeur de psychopédagogie à l'ENNA de Lyon.

[14] François MATRAY, *Pédagogie de l'enseignement technique*, Paris : PUF, 1952, p.91.

### *Du modèle primaire…*

La première phase est celle de l'observation et de la constatation. Il s'agit pour le professeur de faire référence à la matérialité des choses et fonder son enseignement sur l'observation et la manipulation. Dans les instructions officielles comme dans l'ensemble des articles à caractère pédagogique[16] ou dans les leçons-types portant sur les mathématiques[17] publiés dans la revue *TAS*, les élèves doivent « observer des constructions géométriques », ou « la disposition des boutons sur une veste », « décrire », « mesurer », « dénombrer », « réaliser des figures » et « exprimer les résultats des observations ». Ces activités prennent appui sur l'étude de situations concrètes, courantes ou professionnelles, et peuvent parfois mettre en jeu du matériel professionnel que les élèves utilisent à l'atelier comme le palmer ou le calibre à coulisse.

L'enjeu n'est toutefois pas de limiter l'enseignement à la description ou à la constatation d'une propriété, mais d'amener les élèves à comprendre, chercher, critiquer ou contrôler. L'usage des verbes d'actions relevés ici n'est pas anodin. Il témoigne de la volonté des auteurs de rendre l'élève acteur et non plus spectateur dans le processus de transmission des savoirs. Ainsi, le maître doit faire observer les élèves en leur donnant un but et en les faisant opérer avec méthode. Dans un deuxième temps, les élèves doivent réinvestir les résultats obtenus empiriquement sur une multitude d'exemples. Cette deuxième phase vise ainsi à amener les élèves, par un processus de généralisation, à énoncer le principe mathématique expérimenté préalablement. Cette approche intuitive, concrète et expérimentale doit déboucher sur des applications pratiques et contribuer à développer des compétences utiles pour les apprentissages professionnels comme la rigueur, la précision ou le soin pour le dessin technique. En ce sens, les recommandations pédagogiques pour les centres en faveur d'un enseignement mathématique intuitif et inductif, fondé sur l'expérimentation, renvoient à la démarche pédagogique préconisée pour le primaire[18] et mettent en avant l'activité des élèves dans l'élaboration de ses connaissances, pour utiliser un langage actuel, et dans leur réinvestissement dans des applications pratiques.

---

[15] *Ibid.,* p.91.

[16] Fernand CANONGE, *Les centres … op. cit.*; Fernand CANONGE, *Le point de vue … op.cit.* ; Fernand CANONGE, *Les bases … op. cit.* ; Fernand CANONGE, *Les écoles … op. cit. ;* LABARTHE, *L'enseignement… op. cit.*

[17] M ABBES, « Le théorème de Pythagore », *Technique, Art, Science, 2,* 1946, p. 35-37. ; René, CERCELET, « Une leçon de Géométrie dans un Centre d'Apprentissage », *Technique, Art, Science, 1,* 1948, p. 52-54. ; René, CERCELET, « Une leçon de mathématiques dans les Centres d'Apprentissage masculins et féminins », *Technique, Art, Science, 6,* 1949a, p. 52-55. ; René, CERCELET, « Une leçon de mathématiques dans les Centres d'Apprentissage masculins et féminins », *Technique, Art, Science, 8,* 1949b, p. 59-62. ; René, CLUZEL, « Comment préparer et conduire une leçon », *Technique, Art, Science, 1,* 1946, p. 35-37. ; René, CLUZEL, « De la méthode expérimentale à la méthode déductive. *Technique, Art, Science, 2,* 1948, p. 52-54.

### *... au modèle secondaire*

Prendre comme point de départ le concret et l'expérience n'est pas sans rapport avec le fait qu'il s'agit de former les élèves des centres à l'exercice d'un métier technique. Il s'agit, pour un public amené à se heurter à la matérialité des choses dans la pratique professionnelle, de rendre les mathématiques tangibles et utilisables. Toutefois, comme l'indique J. Robert, secrétaire du comité de perfectionnement pédagogique de l'apprentissage en charge de l'écriture des programmes d'enseignement, les enseignements de calcul et de géométrie doivent participer à la formation culturelle des élèves et ne pas être réduits « à des exercices élémentaires purement appliqués au métier, et surtout enseignés suivant des procédés purement empiriques »[19]. En ce sens, les prescriptions insistent sur l'idée que le métier ne doit pas constituer l'unique référence de l'enseignement des mathématiques. D'une part, si pour les élèves des centres, « l'intérêt dominant prend pour objet central le métier »[20], l'enseignement ne doit pas négliger les curiosités désintéressées qui existent chez les apprentis comme chez tous les adolescents. Il doit alors s'ouvrir sur l'histoire, la géographie, le dessin d'art, etc. Par exemple, M. Abbès[21], inspecteur de l'enseignement technique, dans une leçon sur le théorème de Pythagore publiée dans *TAS*, propose des applications se référant à l'histoire de ce théorème. D'autre part, si l'enseignement doit favoriser l'acquisition de connaissances pratiques, sa valeur culturelle réside aussi dans sa participation « à la formation du jugement et du développement de l'intelligence »[22] en conduisant les élèves jusqu'à l'abstraction. L'enchaînement des deux phases observer-généraliser, en somme la démarche inductive, a alors pour objectif de préparer le passage à la troisième phase, la phase déductive. Pour les inspecteurs et les professeurs de centres et d'ENNA, c'est en dépassant les simples constatations lors des démonstrations que l'enseignement de mathématiques peut être à la fois utilitaire et culturel. Par exemple, le programme de géométrie, orienté la première année sur l'étude des figures et la constatation expérimentale des propriétés, vise en deuxième et troisième années à donner à « l'enseignement de la géométrie une orientation et une ampleur telles que le raisonnement y ait une place »[23].

---

[18] Renaud D'ENFERT, « L'enseignement mathématique à l'école primaire de la Troisième république aux années 1960 : enjeux sociaux et culturels d'une scolarisation "de masse" », *Gazette des mathématiciens*, 108, 2006, p. 67-81.

[19] J ROBERT, « La culture dans les Centres d'Apprentissage », *Apprentissage, 2,1945,* p. 95-98

[20] Fernand, CANONGE, *Le point de vue ... op.cit.* p.7

Dans les recommandations pédagogiques, le concret ne s'oppose alors pas à l'abstrait, et n'est pas limité à la manipulation, l'observation ou l'étude de situations pratiques. Un des enjeux de l'enseignement de mathématiques est en effet de donner aux élèves la formation de l'esprit nécessaire pour devenir un travailleur et un citoyen éclairé. Le concret est alors convoqué pour ce qui pose problème et intéresse le public des centres et en vue d'amener à une solution mathématique rationnelle, une explication scientifique : « Le raisonnement doit permettre de dépasser les résultats donnés par l'expérience, celle-ci n'étant utilisée que pour soutenir l'attention, susciter la curiosité »[24]. En ce sens, l'enseignant doit organiser « une classe où l'on fait des mathématiques, où les élèves travaillent sous la direction d'un conseiller qui après avoir bien préparé son travail, conseille et guide »[25] les élèves dans l'élaboration d'une démonstration ou la résolution de problèmes. C'est la méthode interrogative dans le cadre d'un cours dialogué qui est préconisée. L'objectif est de faire participer l'élève par le bais d'un questionnement qui le met dans une position de recherche et qui doit favoriser ainsi la formation de l'esprit.

---

[21] M ABBES, *Le théorème de …op. cit.*

[22] DIRECTION DE L'ENSEIGNEMENT TECHNIQUE, *Instructions sur les programmes et les méthodes des centres d'apprentissage de garçons,* 1945, Lyon : France-empire, p. 41.

[23] *Ibid*, p.42.

[24] LABARTHE, *L'enseignement …op. cit.* p. 45.

[25] FOUSSAT, Réflexions sur l'enseignement des mathématiques et sur les sujets d'examen. *Technique, Art, Science, 4,* 1947, p. 48, p. 48.

Cette démarche pédagogique s'apparente alors à celle pour l'enseignement de mathématiques dans le premier cycle du secondaire[26]. Les élèves seront alors d'autant plus volontiers acteurs que l'enseignant saura susciter leur intérêt :

> *Il y a méthode active lorsque le travail, la recherche, répondent à un besoin véritable chez l'enfant et lorsque l'intérêt, le besoin naissent du dedans et le poussent à conquérir le savoir[27].*

### *L'inscription dans une doctrine pédagogique propre à l'ET*

Fernand Canonge rappelle en 1949 que, dans les centres, « l'apprentissage proprement dit du métier ne sera pas un dressage [...], il devra faire une place importante à l'observation, à l'expérimentation de l'élève, à sa réflexion, à son initiative, à son intelligence de compréhension, de critique, de découverte »[28]. Ainsi, la formation professionnelle doit conduire les élèves à posséder une certaine intelligence des phénomènes technologiques. La pédagogie des mathématiques, fondée sur l'articulation entre les trois phases décrites précédemment (observation, expérimentation, raisonnement) s'apparente alors à la démarche d'apprentissage recommandée dans les revues pour les disciplines professionnelles : « On peut lui appliquer le principe bien connu : observer, comparer, raisonner, induire, déduire et passer ainsi aux applications concrètes »[29].

Certains évoquent ainsi la possibilité de transférer aux mathématiques la démarche pédagogique utilisée dans les disciplines professionnelles et techniques. Ils y voient l'opportunité, moyennant quelques adaptations, de faire profiter cet enseignement de méthodes éprouvées par ailleurs et qui ont montré leur efficacité. Pour d'autres, il s'agit d'établir une sorte d'uniformité des méthodes afin de faciliter le travail des professeurs du domaine professionnel : « La méthode expérimentale est la seule qui permet de bâtir un cours de mathématiques sur lequel peuvent s'appuyer les professeurs de dessin et les maîtres d'atelier »[30].

---

[26] Renaud, D'ENFERT, « Mathématiques modernes et méthodes actives : les ambitions réformatrices des professeurs de mathématiques du secondaire sous la Quatrième République », dans R. D'ENFERT, P. KAHN (dir.), *En attendant la réforme. Disciplines scolaires et politiques éducatives sous la IVe République*, Grenoble, PUG, 2010, p. 115-129.

[27] Roger GAL, Les méthodes actives. *Technique, Art, Science, 7,* 1947, p. 27-28, p. 28.

[28] Fernand CANONGE, *Les écoles … op. cit*. p.41.

[29] L BARON, « L'enseignement professionnel dans les collèges techniques doit être lié à l'organisation industrielle moderne », *Technique, Art, science, 6*, 1949, p.74-76, p.75.

Lors de la création des centres, les méthodes de l'enseignement des mathématiques ne se définissent donc pas simplement par rapport au profil des élèves, mais s'inscrivent dans la doctrine pédagogique de l'ET qui vise à associer action et réflexion. Dans la période d'après-guerre, la définition des principes pédagogiques de l'ET a été fortement influencée par le concept d'humanisme technique[31]. En opposition aux humanités classiques, les promoteurs de l'humanisme technique affirment la valeur culturelle de l'apprentissage du métier. Ils avancent que l'éducation dans les écoles techniques, doit et peut participer à la transmission d'une culture, une culture technique, nécessaire à l'employé qualifié pour s'intégrer dans la société et dépasser sa condition de travailleur à condition que les enseignements et les méthodes soient adaptés à ses possibilités, son caractère et son devenir :

Nos conceptions pédagogiques particulières, tant en ce qui concerne le métier que l'enseignement général, viennent de ce que nous voulons former en même temps, en chacun de nos élèves l'homme et le travailleur.[32]

La culture technique est fondée sur une articulation étroite entre les différents enseignements qui prennent alors pour référence centrale le métier et son apprentissage. Pour Patrice Pelpel et Vincent Troger, la culture technique, telle qu'elle est définie à l'époque, diffère de la culture générale ou scientifique en cela qu'elle est indissociable de l'exercice de l'activité qui la produit[33]. Ses promoteurs souhaitent substituer aux méthodes qu'ils qualifient de traditionnelles, reposant sur la répétition et le cours magistral, et qui visent à « penser pour parler », une pédagogie qui permettrait une interaction entre le savoir et sa mise en pratique, et qui amènerait les élèves à « penser pour agir ». Ils fondent leur réflexion

---

[30] René CLUZEL, *De la méthode … op. cit.* p. 52.
[31] Patrice PELPEL et Vincent TROGER, *Histoire … op. cit.* p. 258.
[32] Fernand CANONGE, *Les écoles … op. cit.* p. 39.
[33] Patrice PELPEL et Vincent TROGER, *Histoire … op. cit.* p. 249-264.

Enseignement des mathématiques, pédagogie concrète et méthodes actives
dans les centres d'apprentissage (1945-19453)

sur la nature pratique et technique de la formation, mais aussi sur l'intérêt « naturel » des jeunes qui fréquentent ces établissements pour le métier et l'influence que ce dernier a sur leur personnalité. Comme le souligne Jean Lamoure, « l'inculcation de cette culture ne pouvant se faire par des méthodes traditionnelles, elle empruntera donc aux méthodes actives »[34].

Les liens qu'entretiennent les cadres et certains acteurs de l'ET comme Paul le Rolland, premier directeur de l'Enseignement technique d'après-guerre, Roger Gal, secrétaire général en 1947 du Groupe français d'éducation nouvelle (GFEN) ou Fernand Canonge, avec le mouvement des pédagogies nouvelles, contribuent à fonder une doctrine pédagogique pour l'enseignement technique fondée sur des « méthodes nouvelles, des méthodes actives, basées sur l'effort, donc sur l'intérêt »[35].

**Un glissement de sens**

Dans le cadre d'une pédagogie qui vise au travers de l'activité à la formation d'un travailleur qualifié et d'un citoyen, et qui est centrée sur la connaissance du profil et du devenir des élèves, Roger Gal rappelle en 1946 que, « dans l'élaboration de nos méthodes, il faudra tenir compte de ces conditions d'âge et de recrutement, en particulier pour nos centres d'apprentissage »[36].

Il ne s'agit donc pas uniquement de définir les méthodes pour l'enseignement de mathématiques autour de ce qui différencie les élèves des écoles techniques des autres élèves du même âge fréquentant le primaire ou le secondaire, ni de les inscrire dans la doctrine pédagogique de l'ET. Il s'agit aussi et surtout de prendre en compte les caractéristiques spécifiques du public des centres, en somme ce qui le distingue de celui des autres établissements techniques. Pour Paul Le Rolland, la pédagogie pour les centres est alors commandée par une donnée générale :

> *Leurs aptitudes manuelles sont plus grandes que leur capacité d'abstraction. Ils sont naturellement plus portés vers les réalisations*

---

[34] Jean LAMOURE, « La revue Technique, Art, Science, 1944-1955 : entre pédagogie et disciplines », dans R. D'ENFERT et P. KAHN (dir.), *En attendant la réforme. Disciplines scolaires et politiques éducatives sous la IVe République* (pp.159-168), Grenoble : PUG, 2010, p. 159-168.

[35] François GOBLOT, cité par Jean LAMOURE, « La revue… », *op. cit.*, p. 162.

[36] Roger GAL, « Notre tâche », *Technique, Art, Science, 1,* 1946, p. 5-6, citation p. 6.

> *concrètes que vers les spéculations intellectuelles.* [De ceci] *découle une certaine méthode pédagogique.*[37]

Cette citation illustre un glissement dans certains discours pédagogiques. Ils ne se réfèrent plus, par exemple, aux centres d'intérêts des élèves ou à la nature professionnelle de la formation, mais à leurs « aptitudes ». Si, dans les centres, les méthodes relatives à l'enseignement des mathématiques doivent s'inscrire dans les lignes directrices évoquées ci-dessus, c'est alors essentiellement autour du profil scolaire, social et psychologique particulier de ces élèves que doit s'élaborer et s'opérationnaliser cette « *certaine méthode pédagogique* ».

### *Convaincre les rétifs : agir pour aimer les mathématiques*

Inspecteurs, professeurs d'ENNA et de centres perçoivent le public des centres comme étant beaucoup plus rétif à l'école que ne l'est celui des autres établissements techniques, et en plus grande difficulté en mathématiques. Par exemple, A. Latelier, enseignant de mathématique, met en avant leur « indifférence ou l'hostilité pour tout ce qui leur rappelle la période scolaire »[38]. Contrairement à leurs camarades des autres établissements techniques, les élèves des centres et leurs familles, souvent d'origine très modeste, accordent moins d'importance aux études et ne voient dans les centres que la possibilité d'apprendre rapidement un métier. Les problèmes comportementaux que rencontrent les enseignants[39] ainsi que les faibles résultats dans les enseignements généraux semblent alors indiquer, pour les acteurs proches du terrain, que les élèves ne sont attirés que par la formation professionnelle et qu'ils délaissent les autres disciplines qu'ils jugent parfois inutiles[40]. La nature professionnelle de la formation et son organisation, principalement axée sur la pratique qui occupe plus de la moitié de l'emploi du temps, ainsi que le poids des disciplines professionnelles à l'examen, y contribuent.

---

[37] Paul LE ROLLAND cité par AVENIR, « Les carrières de l'Enseignement Technique », *Avenir, 103-104*, 1959, p. 137-154, p. 137.

[38] A LATELIER, « Arithmétique, Amortissement d'une dette d'un matériel, d'un outillage », *Apprentissage, 5*, 1947, p.255-256, p. 256.

[39] Patrice PELPEL et Vinent TROGER, *Histoire… op. cit.* p.170-171.

[40] René, BATAILLON, « Le but et l'esprit des enseignements scientifiques », dans Ministère de l'Éducation nationale, *Encyclopédie générale de l'éducation française*, Paris, Rombaldi, 1954, p. 250 – 255.

Enseignement des mathématiques, pédagogie concrète et méthodes actives
dans les centres d'apprentissage (1945-19453)

Si le recours aux méthodes actives centrées sur le concret s'inscrit toujours dans le cadre d'une rupture avec les méthodes traditionnelles, l'enjeu est moins de se référer à la nature professionnelle de la formation que de convaincre les élèves de l'utilité de faire des mathématiques, de ne pas les ennuyer ou les dégoûter du livre et de la connaissance gratuite.

*Il convient, pour intéresser les jeunes apprentis et fixer leur esprit, de rattacher autant que possible l'objet du cours à une notion professionnelle*[41]

L'activité des élèves se réfère à l'utilisation qu'un professionnel ou un citoyen peut faire des mathématiques. Observer et expérimenter visent à montrer aux élèves comment on se sert des mathématiques pour la résolution de problèmes concrets, qu'ils soient courants, professionnels ou d'examens, et non plus dans la perspective d'un passage vers l'abstraction.

Ainsi, pour René Cercelet, inspecteur de l'ET, il faut que les élèves « sentent par les applications nombreuses que nous leur donnons (constructions graphiques, étude d'appareils divers, problèmes professionnels) que les mathématiques sont un instrument de travail au même titre que la lime, le rabot ou la perceuse »[42]. L'enjeu est aussi de donner du sens à l'enseignement en établissant une correspondance entre les différentes leçons du cours de mathématiques et leur utilité dans le domaine pratique afin de permettre aux élèves « de VOIR et de voir exactement »[43] les applications pratiques de la notion mathématique étudiée.

L'utilisation de l'observation et de l'expérimentation, dans le cadre d'une démarche inductive, est aussi perçue comme le moyen de s'affranchir des prérequis mathématiques. En effet, les premières instructions de 1945 et les recommandations des revues[44] laissent aux professeurs la possibilité de déroger à la progression « normale » de leur enseignement en leur laissant une grande liberté dans l'agencement des leçons, et en les engageant à utiliser une démarche inductive pour donner « parfois les connaissances nécessaires sous une forme simplifiée, intuitive, provisoire, quitte à y revenir plus tard pour les préciser et les compléter »[45]. Il s'agit de mettre en adéquation cet enseignement avec les besoins des disciplines professionnelles afin d'apporter aux é-lèves l'appui des connaissances mathématiques nécessaires au moment où ils en ont besoin.

---

[41] René CERCELET, *Une leçon…, op. cit* p. 54.

[42] *Ibid* p.52.

[43] Henri GIDROL, Un essai d'enseignement concret des mathématiques au Centre d'Apprentissage. *Technique, Art, Science, 9,* 1955, p. 31-37, p, 33.

[44] René ORIOL, « Réflexion sur l'enseignement des mathématiques dans les centres d'apprentissages », *Technique, Art, Science, 1*, 1949, p. 5-7.

### *Un enseignement formateur de l'esprit pour tous, même pour les concrets*

Dans les revues, les qualificatifs pour désigner les élèves des centres ne manquent pas. Ils sont souvent perçus comme « incapables d'abstraction »[46], possédant une intelligence « concrète »[47], et « d'âge intellectuel très faible ou encore retardés »[48]. Les praticiens soulignent la difficulté d'opérationnaliser un enseignement de mathématiques formateur de l'esprit pour des classes hétérogènes composées en majorité d'« élèves plus ou moins retardés qui constituent le gros de leur effectif »[49]. Afin d'étudier la prise en charge dans les recommandations pédagogiques de ces élèves perçus comme peu doués pour l'abstraction et « qui ne peuvent pas apprendre »[50], nous nous sommes intéressé aux parties « démonstration d'une propriété »[51] et « passage à l'abstraction »[52] des leçons intitulées *Une leçon de géométrie dans un centre d'apprentissage* et *Le théorème de Pythagore*, publiées respectivement par René Cercelet et M. Abbès. Nous rappelons que cette phase de la démarche pédagogique préconisée pour l'enseignement des mathématiques au travers d'un cours dialogué. Le moyen utilisé par les professeurs pour amener les élèves à généraliser le résultat obtenu à partir des premières constatations consiste à multiplier les exemples et à les guider au moyen d'un questionnement dont nous avons relevé quelques extraits :

*Ne voyez-vous pas dans le triangle que vous avez construit un ou des angles plats ? Infailliblement, les élèves en montreront une infinité.*

---

[45] René Bataillon, *Le but et… op. cit.* p. 251.

[46] Direction de l'enseignement technique, *Instructions …, op. cit.* p.13.

[47] René, Cercelet, *Une leçon de Géométrie… op.cit.* p.52.

[48] *Ibid*, p. 52.

[49] Jean Alphand, « Arithmétique », *Apprentissage* 9 (1946), p.726.

[50] Labarthe, *L'enseignement …op. cit.* p. 45.

[51] M Abbes, *Le théorème… op. cit.,* p. 36.

[52] René Cercelet, *Une leçon…, (1949) op. cit.,* p. 54.

Enseignement des mathématiques, pédagogie concrète et méthodes actives
dans les centres d'apprentissage (1945-19453)

*Les élèves suivent alors les explications du professeur au tableau, qui doit pouvoir faire formuler par ses élèves la phrase capitale suivante : …*

 *Le professeur doit alors faire ressortir l'admirable logique du raisonnement et la rigueur de cette démonstration basée sur une construction simple que l'on est amené à réaliser après la constatation de la propriété. […] Le professeur doit faire ressortir le merveilleux de cette constatation.*[53]

*Avec aisance les élèves généralisent* [54]

*La comparaison vous conduit, même avec les plus médiocres à la formule…* [55]

Nous pouvons constater que l'activité de recherche des élèves semble très réduite voire inexistante dans les exemples proposés. En effet, dans leur formulation, les questions posées contiennent les réponses que les élèves doivent trouver. Le dialogue entre le professeur et les élèves possède ainsi un aspect artificiel voire idéal. Le problème des élèves « médiocres », « peu doués pour l'abstraction », souligné à plusieurs reprises dans l'ensemble des discours, semble ici se résoudre de lui-même grâce à l'émerveillement que ne doit pas manquer de produire la démarche inductive chez les élèves. Il s'agit alors pour l'enseignant moins de les laisser chercher et participer à l'élaboration du savoir que de les guider à l'aide d'un questionnement judicieux afin qu'ils ne restent jamais trop longtemps seuls face à une difficulté en leur donnant l'impression que ce sont eux qui résolvent les problèmes[56]. Il s'agit alors, pour Labarthe, chargé de cours à l'ENNA, de « ne pas tout démontrer systématiquement, mais seulement lorsque l'élève exige la démonstration ou est en mesure de la comprendre »[57]. L'activité de l'élève est ainsi circonscrite à la manipulation et l'expérimentation et ne s'étend pas jusqu'au raisonnement mathématique, limitant ainsi l'enseignement au concret pour des « esprits concrets ».

C'est donc un modèle pédagogique original qui est proposé pour l'enseignement de mathématiques dans l'enseignement technique court. En référence aux caractéristiques, sociales, psychologiques, scolaires particulières de ce nouveau public scolaire à qui ne convient pas les méthodes « traditionnelles », les propositions s'inscrivent dans le courant des pédagogies nouvelles et visent à promouvoir un enseignement actif, concret et qui prenne en compte les centres d'intérêts des élèves. Les termes « actif » et « concret » renvoient à des acceptions différentes selon que les discours font références aux facettes utilitaire ou désintéressée de l'enseignement. Dans une dialectique revendiquée à l'époque par l'ET, entre action et réflexion, la démarche pédagogique pour l'enseignement des mathématiques s'inscrit dans le passage, sinon dans l'articulation, entre le modèle primaire qui vise à l'acquisition d'une culture pratique et utilitaire, et celui du premier cycle du secondaire qui doit conduire à l'abstraction.

---

[53] René CERCELET, *Une leçon…, (1949) op. cit.,* p. 54.

[54] *Ibid*, p. 52.

[55] M ABBES, *Le théorème… op. cit.,* p. 36.

[56] P BETRAMA, « Comment renouveler l'intérêt d'un problème déjà étudié et mal connu », *Technique, Art, Science, 4,* 1947, p.52-54.

[57] LABARTHE, *L'enseignement …op. cit.* p. 45

En référence au profil des élèves des centres et en particulier à leur supposée incapacité à accéder à l'abstraction, un glissement de sens s'opère dans les discours tenus par des acteurs proches du terrain et en lien avec l'opérationnalisation de l'enseignement. Le concret ne renvoie plus à une démarche pédagogique mais à la tournure d'esprit des élèves et, pour reprendre une citation de Roger Gal, les « méthodes actives » sont devenues les « méthodes vivantes »[58] renvoyant à l'action physique des élèves, loin de toute élaboration de connaissances.

À la fois l'objet et les enjeux de ce travail contribuent à étendre et réexplorer les questionnements présents dans cet ouvrage sur l'enseignement de mathématiques dans le primaire dans la période d'après-guerre, et notamment l'activité des élèves ou leur rapport au concret. Il contribue ainsi à combler un vide dans l'historiographie de l'enseignement mathématiques dans la seconde moitié du XXe siècle, principalement centrée sur le primaire et le secondaire, en interrogeant notamment l'articulation entre primaire, secondaire et technique court.

---

[58] Roger GAL, *Les méthodes … op. cit.* p. 28.

Enseignement des mathématiques, pédagogie concrète et méthodes actives dans les centres d'apprentissage (1945-19453)

Index des noms

Abbès, M ; Alphand, Jean ; Baron, L. ; Bataillon, René ; Besse, Laurent ; Bétrama, P. ; Bourjade, Jean ; Brucy, Guy ; Canonge, Fernand ; Cercelet, René ; Cluzel, René ; d'Enfert, Renaud ; Dintzer, Louis ; Foussat ; Pierre Kahn; Gal, Roger ; Gidrol, Henri ; Goblot François, Gutierrez, Laurent ; Juif, Paul ; Labarthe ; Lamoure, Jean ; Latelier, A ; Lebeaume, Joël ; Le Rolland, Paul ; Lopez, Maryse ; Maillard, Fabienne ; Matray, François ; Moreau, Oriol, R ; Gilles ; Pelpel, Patrice ; Prost, Antoine ; Robert, A ; Sido, Xavier ; Troger, Vincent.